\documentclass[a4paper]{amsart}
\usepackage[ps2pdf]{hyperref}
\usepackage{amsmath,amssymb}
\usepackage{bbold}
\usepackage{psfrag}
\usepackage{graphicx}
\usepackage[all,poly]{xy}
\usepackage[latin1]{inputenc}

\newtheorem{thm}{Theorem}[section]

\newcommand{\ttt}{\mathcal{T}}
\newcommand{\WRT}{\mathrm{WRT}}
\newcommand{\co}{\colon}
\newcommand{\id}{\mathrm{id}}
\newcommand{\cc}{\mathcal{C}}
\newcommand{\bb}{\mathcal{B}}
\newcommand{\dd}{\mathcal{D}}
\newcommand{\zz}{\mathcal{Z}}
\newcommand{\CC}{\mathbb{C}}
\newcommand{\RR}{\mathbb{R}}
\newcommand{\kk}{\Bbbk}
\newcommand{\trait}{\nobreakdash-\hspace{0pt}}
\newcommand{\Z}{\mathbb{Z}}
\newcommand{\un}{\mathbb{1}}
\newcommand{\End}{\mathrm{End}}
\newcommand{\Hom}{\mathrm{Hom}}
\newcommand{\Vect}{{\mathrm{Vect}}}
\newcommand{\tr}{\mathrm{tr}}
\newcommand{\lev}{\mathrm{ev}}
\newcommand{\rev}{\widetilde{\mathrm{ev}}}
\newcommand{\lcoev}{\mathrm{coev}}
\newcommand{\rcoev}{\widetilde{\mathrm{coev}}}
\newcommand{\rsdraw}[3]{\raisebox{-#1\height}{\scalebox{#2}{\includegraphics{#3.eps}}}}

\begin{document}

\title{3d TQFTs and 3-manifold invariants}
\author{K\"{u}r\c{s}at S\"{o}zer}
 \address{K\"{u}r\c{s}at S\"{o}zer\newline
  \indent  McMaster University, Hamilton, ON L8S 4E8, Canada\\}
\email{sozerk@mcmaster.ca}
\author{Alexis Virelizier}
\address{Alexis Virelizier\newline
\indent Univ. Lille, CNRS, UMR 8524 - Laboratoire Paul Painlev\'e, F-59000 Lille, France\\}
\email{alexis.virelizier@univ-lille.fr}

\date{\today}

\begin{abstract}
This is an invited contribution to the 2nd edition of the Encyclopedia of Mathematical Physics.
We give an overview of 3-dimensional topological quantum field theories (TQFTs) and the corresponding quantum invariants of 3-manifolds. We recall the main algebraic concepts and constructions, such as modular and spherical fusion categories, the Witten-Reshetikhin-Turaev and Turaev-Viro theories, and the relation between these two TQFTs. We also briefly discuss generalizations of these constructions by providing a (non-exhaustive) review of some recent works on 3-dimensional extended TQFTs, defect TQFTs, homotopy QFTs, and non-semisimple TQFTs.
\end{abstract}

\maketitle

\setcounter{tocdepth}{1}
\tableofcontents

\section{Introduction}\label{sec-Intro}
Topological invariants are quantities associated with a topological space that do not change under continuous deformations of the space. One way to determine if two spaces are topologically distinct from each other is to compare the values of these invariants. Many topological invariants (such as (co)homology and homotopy theories) have been introduced and thoroughly studied since the XIX$^{\mathrm{th}}$ century, allowing for a complete classification of several families of topological objects. However, the study of 3-dimensional manifolds using topological invariants remains  a very active area of research. In particular, a new class of topological invariants of 3-manifolds, called {\it quantum invariants}, emerged in the 1980s.

Quantum invariants originate from the idea of relating the topology of smooth manifolds to the partition functions of certain quantum field theories (QFTs). This idea was first proposed by Schwarz in 1978, and elaborated by Witten in 1988 who showed that the Chern-Simons QFT can produce the Jones polynomial, a polynomial invariant of knots and links in the $3$-sphere. This was the beginning of a fascinating interaction between mathematics and theoretical physics. Witten also conjectured that Chern-Simons theory can be used to define more general invariants of $3$-manifolds, which were later constructed rigorously by Reshetikhin and Turaev in 1989 using quantum groups. These quantum invariants, more generally defined using modular categories, are known as the Witten-Reshetikhin-Turaev invariants and extend to $3$\trait dimensional topological quantum field theories (TQFTs), which are QFTs that depend  only on the topology and not on the geometry of the manifolds.
A second important family of quantum $3$-manifold invariants comes from the Turaev-Viro-Barrett-Westbury  state sum construction on triangulations of $3$-manifolds, defined in  the  1990s using fusion categories. The quantum field theory motivating these state sum invariants is the Ponzano-Regge model for $3$-dimensional lattice gravity.
Since then, quantum $3$\trait manifold invariants and their associated TQFTs have been extensively studied and successfully generalized in several directions (including extended TQFTs, defect TQFTs, homotopy QFTs, non-semisimple TQFTs).

This review is organized as follows. Section~\ref{sect-categories} is  dedicated to algebraic preliminaries on the pivotal, ribbon, fusion, and modular categories. Section~\ref{sect-TQFTs} is devoted to generalities on 3-dimensional TQFTs. In Section~\ref{sect-surgery-TQFT}, we define  the Witten-Reshetikhin-Turaev surgery invariants and their TQFTs from modular categories. In Section~\ref{sect-state-sum-TQFT}, we define the state sum invariants and their TQFTs from spherical fusion categories. Section~\ref{sect-comparison} is devoted to the comparison of surgery and state sum approaches. Finally, in Section~\ref{sect-generalizations}, we end with a (non-exhaustive)  review of some more recent works on extended TQFTs, defect TQFTs, homotopy QFTs, and  non-semisimple TQFTs.

\section{Algebraic preliminaries}\label{sect-categories}
The definition of quantum invariants begins with fixing suitable algebraic data, which are best described in terms of monoidal categories.

\subsection{Pivotal categories}\label{sect-pivotal}
Let $\cc$ be a monoidal category (i.e., a category with an associative tensor product and a unit object $\un$).
A \emph{left duality} in $\cc$ assigns to any object $X$ of $\cc$ an object~$X^*$ together with two morphisms $\lev_X \co X^* \otimes  X \to \un$ and $\lcoev_X \co \un \to X \otimes X^*$ in $\cc$ (the left evaluation and coevaluation) such that
\begin{align*}
& (\id_X \otimes \lev_X)(\lcoev_X \otimes \id_X)=\id_X, \\
& (\lev_X \otimes \id_{X^*})(\id_{X^*} \otimes \lcoev_X)=\id_{X^*}.
\end{align*}
The left dual of a morphism $f \co X \to Y$ is then the morphism $f^*\co Y^* \to X^*$ defined by
$$
f^*=(\lev_Y \otimes \id_{X^*})(\id_{Y^*} \otimes f \otimes \id_{X^*})(\id_{Y^*} \otimes \lcoev_X).
$$
We will often abstain (by abuse) from writing down the following canonical isomorphisms
\begin{align*}
&X^{**}\cong X, && (X \otimes Y)^* \cong Y^* \otimes X ^*, && \un^* \cong \un.
\end{align*}

A \emph{pivotal structure} in $\cc$ is a left duality in $\cc$ together with a natural isomorphism $\phi=\{\phi_X \co X \to X^{**}\}_{X \in \cc}$ which is monoidal in the sense that $\phi_{X \otimes Y}=\phi_X \otimes \phi_Y$. The right evaluation and coevaluation associated with an object $X \in \cc$ are then defined by
\begin{align*}
&\rev_X=\lev_{X^*}(\phi_X \otimes \id_{X^*}) \co X \otimes X^* \to \un,\\
&\rcoev_X=(\id_{X^*} \otimes \phi^{-1}_X)\lcoev_{X^*}\co \un \to X^* \otimes X.
\end{align*}
The (co)evaluation morphisms allow to define the \emph{left trace} and \emph{right trace} of any endomorphism $g \co X \to X$ as
\begin{align*}
&\tr_l(g)=\lev_X(\id_{X^*} \otimes g) \rcoev_X \co \un \to \un,\\
&\tr_r(g)=\rev_X(g \otimes \id_{X^*}) \lcoev_X \co \un \to \un.
\end{align*}
Both take values in the commutative monoid $\End_\cc(\un)$ of endomorphisms of the monoidal unit $\un$ and share a number of properties of the standard trace of matrices such as
$\tr_l(fh)=\tr_l(hf)$ and $\tr_l(g)=\tr_r(g^*)=\tr_l(g^{**})$ (and similarly  with $l,r$ exchanged).
The \emph{left} and \emph{right dimensions} of an object   $X\in\cc$
are defined by
$$
\dim_l(X)=\tr_l(\id_X) \quad \text{and} \quad \dim_r(X)=\tr_r(\id_X).
$$
Note that isomorphic objects have the same
dimensions and $\dim_l(\un)=\dim_r(\un)=\id_{\un}$.

A \emph{pivotal category} is a monoidal category endowed with a pivotal structure.

\subsection{Penrose graphical calculus}\label{sect-Penrose}
We represent morphisms in a pivotal category~$\cc$ by planar  diagrams to be read from the bottom to the top.
Diagrams are made of   oriented arcs colored by objects of~$\cc$  and of boxes colored by morphisms of~$\cc$.  The arcs connect
the boxes and   have no mutual intersections or self-intersections.
The identity $\id_X$ of an object $X$, a morphism $f\co X \to Y$, the composition of two morphisms $f\co X \to Y$ and $g\co Y \to
Z$, and the monoidal product of two morphisms $\alpha \co X \to Y$
and $\beta \co U \to V$  are represented as follows:
\begin{center}
\psfrag{X}[Bc][Bc]{\scalebox{.7}{$X$}} \psfrag{Y}[Bc][Bc]{\scalebox{.7}{$Y$}} \psfrag{h}[Bc][Bc]{\scalebox{.8}{$f$}} \psfrag{g}[Bc][Bc]{\scalebox{.8}{$g$}}
\psfrag{Z}[Bc][Bc]{\scalebox{.7}{$Z$}} $\id_X=$ \rsdraw{.45}{.9}{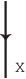}\,,\quad $f=$ \rsdraw{.45}{.9}{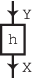} ,\quad $gf=$ \rsdraw{.45}{.9}{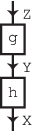}\,, \quad
\psfrag{X}[Bc][Bc]{\scalebox{.7}{$X$}} \psfrag{h}[Bc][Bc]{\scalebox{.8}{$\alpha$}}
\psfrag{Y}[Bc][Bc]{\scalebox{.7}{$Y$}}  $\alpha\otimes \beta=$ \rsdraw{.45}{.9}{morphism} \psfrag{X}[Bc][Bc]{\scalebox{.8}{$U$}} \psfrag{h}[Bc][Bc]{\scalebox{.8}{$\beta$}}
\psfrag{Y}[Bc][Bc]{\scalebox{.7}{$V$}} \rsdraw{.45}{.9}{morphism}\,.
\end{center}
A box whose lower/upper side has no attached strands represents a morphism with source/target $\un$.
If an arc colored by $X$ is oriented upward,
then the corresponding object   in the source/target of  morphisms
is $X^*$. For example, $\id_{X^*}$  and a morphism $f\co X^* \otimes
Y \to U \otimes V^* \otimes W$  may be depicted as:
\begin{center}
 $\id_{X^*}=$ \, \psfrag{X}[Bl][Bl]{\scalebox{.7}{$X$}}
\rsdraw{.45}{.9}{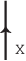} $=$  \,
\psfrag{X}[Bl][Bl]{\scalebox{.7}{$X^*$}}
\rsdraw{.45}{.9}{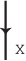}  \quad and \quad
\psfrag{X}[Bc][Bc]{\scalebox{.7}{$X$}}
\psfrag{h}[Bc][Bc]{\scalebox{.8}{$f$}}
\psfrag{Y}[Bc][Bc]{\scalebox{.7}{$Y$}}
\psfrag{U}[Bc][Bc]{\scalebox{.7}{$U$}}
\psfrag{V}[Bc][Bc]{\scalebox{.7}{$V$}}
\psfrag{W}[Bc][Bc]{\scalebox{.7}{$W$}} $f=$
\rsdraw{.45}{.9}{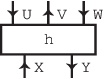} \,.
\end{center}
The duality morphisms are depicted as
\begin{center}
\psfrag{X}[Bc][Bc]{\scalebox{.7}{$X$}} $\lev_X=$ \rsdraw{.45}{.9}{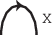}\,,\quad
 $\lcoev_X=$ \rsdraw{.45}{.9}{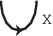}\,,\\
$\rev_X=$ \rsdraw{.45}{.9}{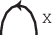}\,,\quad
\psfrag{C}[Bc][Bc]{\scalebox{.7}{$X$}} $\rcoev_X=$
\rsdraw{.45}{.9}{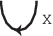}\,.
\end{center}
The dual of a morphism $f\co X \to Y$ can be depicted as
\begin{center}
\psfrag{X}[Bc][Bc]{\scalebox{.7}{$X$}} \psfrag{h}[Bc][Bc]{\scalebox{.8}{$f$}}
\psfrag{Y}[Bc][Bc]{\scalebox{.7}{$Y$}} \psfrag{g}[Bc][Bc]{\scalebox{.8}{$g$}}
$f^*=$ \rsdraw{.45}{.9}{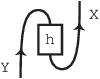}$=$ \rsdraw{.45}{.9}{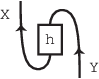}
\end{center}
and the traces of an endomorphism $g\co X \to X$ as
\begin{center}
\psfrag{X}[Bc][Bc]{\scalebox{.7}{$X$}}  \psfrag{g}[Bc][Bc]{\scalebox{.8}{$g$}}
$\tr_l(g)=$ \rsdraw{.45}{.9}{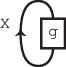}\,,\quad  $\tr_r(g)=$ \rsdraw{.45}{.9}{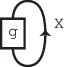}\,.
\end{center}
Note that the morphisms represented by the diagrams
are invariant under isotopies of the diagrams in the plane keeping
fixed the bottom and top endpoints (see \cite{JS,TVi5}).

\subsection{Spherical categories}\label{sect-spherical}
A \emph{spherical category} is a pivotal category whose left and
right traces are equal, i.e.,  $\tr_l(g)=\tr_r(g)$ for every
endomorphism $g$ of an object. Then $\tr_l(g)$ and $ \tr_r(g)$ are
denoted $\tr(g)$ and called the \emph{trace of $g$}. In particular, the
left and right dimensions of an object~$X$ are equal, denoted   $\dim(X)$,
and called the \emph{dimension of $X$}.

For spherical categories, the corresponding Penrose gra\-phi\-cal calculus has the following property:   the morphisms represented by
 diagrams are invariant under isotopies of   diagrams in the 2-sphere $S^2=\RR^2\cup
\{\infty\}$, i.e., they are preserved under isotopies pushing   arcs of
the diagrams across~$\infty$.  For example, the diagrams above
representing $\tr_l(g)$ and $\tr_r(g)$ are related by such an
isotopy. Note that the condition $\tr_l(g)=\tr_r(g)$ for all $g$ is therefore
necessary (and in fact sufficient) to ensure this property.

\subsection{Braided categories}\label{sect-braided}
A \emph{braiding} in a monoidal category $\bb$ is a natural isomorphism
$ c=\{c_{X,Y} \co  X \otimes Y\to Y \otimes X\}_{X,Y \in
\bb} $ such that
\begin{align*}
&c_{X, Y\otimes Z}=(\id_Y \otimes c_{X, Z})(c_{X, Y} \otimes \id_Z), \\
&c_{X\otimes Y,Z}=(c_{X,Z} \otimes \id_Y)(\id_X \otimes c_{Y,Z}).
\end{align*}
These conditions imply that $c_{X,\un}=c_{\un,X}=\id_X$ for
any object $X$. A monoidal category endowed with a braiding is said to be \emph{braided}.

Let $\bb$ be a braided pivotal category. The braiding  and its inverse are depicted
as
\begin{center}
\psfrag{X}[Bc][Bc]{\scalebox{.8}{$X$}}
\psfrag{Y}[Bc][Bc]{\scalebox{.8}{$Y$}}
$c_{X,Y}=\,$\rsdraw{.45}{.9}{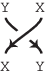} \quad \text{and} \quad
$c^{-1}_{Y,X}=\,$\rsdraw{.45}{.9}{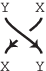}.
\end{center}
The family $\theta=\{\theta_X \co X \to X\}_{X \in \bb}$, defined by
$$ \theta_X  =
\psfrag{X}[Bc][Bc]{\scalebox{.8}{$X$}}\rsdraw{.45}{.9}{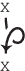}\,=(\id_X
\otimes \rev_X)(c_{X,X} \otimes \id_{X^*})(\id_X \otimes
\lcoev_X),
$$
is called \emph{twist} of $\bb$. It is a natural isomorphism and
satisfies  $\theta_\un=\id_\un$ and
$$\theta_{X\otimes Y}=(\theta_X \otimes \theta_Y)c_{Y,X}c_{X,Y}.$$

\subsection{Ribbon categories}\label{sect-ribbon}
A \emph{ribbon category} is a braided pivotal category $\bb$ whose
twist  is self-dual, i.e., $(\theta_X)^*=\theta_{X^*}$ for all objects $X$ of $\bb$. This condition  is equivalent to the equality of morphisms
$$
\psfrag{X}[Bc][Bc]{\scalebox{.8}{$X$}} \rsdraw{.45}{.9}{theta1}
=\rsdraw{.45}{.9}{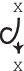}.
$$
The inverse of the twist is then computed by
$$
\theta_X^{-1}=\psfrag{X}[Bc][Bc]{\scalebox{.8}{$X$}} \rsdraw{.45}{.9}{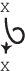}=\rsdraw{.45}{.9}{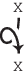}.
$$

Ribbon categories are spherical and nicely fit into the theory of knots and links in $S^3$. A \emph{link} $L \subset S^3$ is a closed one-dimensional submanifold of $S^3$. (A manifold is \emph{closed} if it is compact and has no boundary.) A link is \emph{oriented} (resp.\@ \emph{framed}) if all its components are oriented (resp.\@ provided with a homotopy class of nonsingular normal vector fields). Any ribbon category $\bb$  gives rise to an invariant of $\bb$-colored framed oriented links in $S^3$. Here, a link is \emph{$\bb$-colored} if each of its components is endowed with an object of $\bb$ (called the \emph{color} of this component). Namely,   every   $\bb$-colored framed oriented link   $L\subset S^3$ determines an endomorphism of the unit object
$$\langle L \rangle_\bb\in \End_\bb (\un)$$ which turns out to be an isotopy invariant of $L$. To compute $\langle L \rangle_\bb$, present $L$ by a planar diagram with only double transversal crossings such that the framing of~$L$ is orthogonal to the plane, and then apply the Penrose graphical calculus  to this $\bb$-colored diagram (using the braiding and its inverse for the positive and negative crossings). The axioms of a ribbon category imply that $\langle L \rangle_\bb$ does not depend on the chosen plane diagram for $L$. For example, $$\langle O_X \rangle_\bb = \dim(X)$$ for the trivial knot $O_X$ with zero framing and color $X \in \bb$.

Further constructions need the notion of a tangle. An (oriented) \emph{tangle} is a compact (oriented) one-dimensional submanifold of $\RR^2 \times [0, 1]$ with endpoints on $\RR \times 0 \times \{0, 1\}$. Near each of its endpoints, an oriented tangle $T$ is directed either down or up, and thus acquires a sign $\pm 1$. Then one can view $T$ as a morphism from the sequence of $\pm 1$'s associated with its bottom ends to the sequence of $\pm 1$'s associated with its top ends. Tangles can be composed by putting one on top of the other. This defines a monoidal category of tangles $\ttt$ whose objects are finite sequences of~$\pm 1$'s and whose morphisms are isotopy classes of framed oriented tangles. Given a ribbon category $\bb$, we can consider $\bb$-colored tangles, that is, (framed oriented) tangles whose components are labeled with objects of $\bb$. They form a category $\ttt_{\bb}$. Links appear here as tangles without endpoints, that is, as morphisms $\emptyset \to \emptyset$. The link invariant $\langle L \rangle_\bb$ generalizes to a functor $\langle \cdot \rangle_\bb \co \ttt_\bb \to \bb$, see \cite{Tu1}.

\subsection{Fusion categories}\label{sect-fusion}
Let $\kk$ be a field.  A monoidal category is \emph{$\kk$-linear} if its Hom sets are $\kk$-vector spaces, and the composition and monoidal product of morphisms are $\kk$-bilinear. Such a category is \emph{additive} if any finite family of objects has a direct sum.

An object~$S$ of a $\kk$-linear monoidal category $\cc$ is \emph{simple} if the $\kk$-vector space $\End_\cc(S)$ is one dimensional. Then the map $\kk \to \End_\cc(S)$, $k \mapsto k \, \id_S$  is a $\kk$-algebra isomorphism. It is used to identify $\End_\cc(S)=\kk$.

A \emph{fusion category} (over $\kk$) is an additive $\kk$-linear pivotal category~$\cc$ such that
each object of $\cc$ is a (finite) direct sum of simple objects, $\Hom_\cc(i,j)=0$ for any non-isomorphic simple objects $i,j  $ of $\cc$, the unit object $\un$ is simple, and the set of isomorphism classes of simple objects of~$\cc$   is finite.
These conditions imply that all the $\Hom$ spaces in $\cc$ are finite dimensional $\kk$-vector spaces.

In a fusion category, the left and right dimensions of any simple object of~$\cc$ are nonzero in~$\End_\cc(\un)=\kk$. Also, a fusion category is spherical if and only if any simple object has equal left and right dimensions.

The \emph{dimension} of a fusion category $\cc$ is
\begin{equation*}
\dim (\cc)=\sum_{i \in I} \dim_l(i) \dim_r(i) \in \kk,
\end{equation*}
where $I$ is any representative set of simple objects of  $\cc$ (meaning that $\un \in I$ and every simple object of
$\cc$ is isomorphic to a unique element of~$I$). By \cite{ENO}, $\dim (\cc)\neq 0$ when  $\kk $ is an algebraically closed field of characteristic zero. For spherical~$\cc$, we have $\dim (\cc)=\sum_{i \in I} \dim( i)^2$.

A standard example of a spherical fusion category with nonzero dimension is the   category of finite dimensional representations (over~$\kk$) of a finite group  whose order is relatively prime to  the characteristic of~$\kk$. More interesting examples of spherical fusion categories
are derived from the theory of subfactors, see \cite{EK,KS2}.

\subsection{$6j$-symbols}\label{sect-6j}
The $6j$-symbols were first introduced by the physicists Wigner and Racah in the theory of representations of  $SU_2(\CC)$.  The $6j$-symbols  have been extensively used in the  theory of angular momentum in quantum mechanics and in the Ponzano-Regge approach to  quantum
gravity in dimension three. Also, the $6j$-symbols   play a special role  in  3-dimensional state sum TQFTs (see Section~\ref{sect-state-sum-TQFT}). We will need two versions of the $6j$-symbols (among the $2^6=64$ versions, each of them corresponding  to a choice  of orientation for the edges of a tetrahedron, see \cite{TVi5}).

Let $\cc$ be a spherical fusion category and $I$ be a representative set of simple objects of  $\cc$. For $i,j,k \in I$, consider the multiplicity spaces
$$
H_{i,j}^k=\Hom_\cc(i \otimes j,k) \quad \text{and} \quad H_k^{i,j}=\Hom_\cc(k,i \otimes j).
$$
The positive $6j$-symbol associated with  $i,j,k,\ell, m, n \in I$ is the $\kk$-linear form
$$
\begin{Bmatrix}i&j&k \\ \ell&m&n\end{Bmatrix}_{\!+} \co H_m^{k,\ell}\otimes H_{j,\ell}^n \otimes   H_n^{i,m} \otimes  H_{i,k}^j \to \kk
$$
which maps $\alpha \otimes \beta \otimes \gamma \otimes \delta$ to
$$
\psfrag{a}[Bc][Bc]{\scalebox{.9}{$\alpha$}}
\psfrag{b}[Bc][Bc]{\scalebox{.9}{$\beta$}}
\psfrag{c}[Bc][Bc]{\scalebox{.9}{$\gamma$}}
\psfrag{d}[Bc][Bc]{\scalebox{.9}{$\delta$}}
\psfrag{i}[Br][Br]{\scalebox{.7}{$i$}}
\psfrag{j}[Br][Br]{\scalebox{.7}{$j$}}
\psfrag{k}[Bl][Bl]{\scalebox{.7}{$k$}}
\psfrag{l}[Bl][Bl]{\scalebox{.7}{$\ell$}}
\psfrag{m}[Bl][Bl]{\scalebox{.7}{$m$}}
\psfrag{n}[Bl][Bl]{\scalebox{.7}{$n$}}
\tr\bigl(\beta(\delta \otimes \id_\ell)(\id_i \otimes \alpha)\gamma \bigr) = \, \rsdraw{.45}{.9}{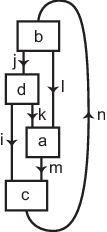} \;.
$$
Similarly, the negative $6j$-symbol is the $\kk$-linear form
$$
\begin{Bmatrix}i&j&k \\ \ell &m&n\end{Bmatrix}_{\!-}  \co H^m_{k,\ell}\otimes H^{j,\ell}_n \otimes   H^n_{i,m} \otimes  H^{i,k}_j \to \kk
$$
defined by $\alpha \otimes \beta \otimes \gamma \otimes \delta \mapsto \tr\bigl(\gamma(\id_i \otimes \alpha)(\delta \otimes \id_\ell)\beta \bigr)$.

Note that if the multiplicity spaces are at most one dimensional and have canonical basis elements (as in the $SU_2(\CC)$ case), then the  $6j$-symbols can be interpreted as numbers.

The $6j$-symbols   satisfy  beautiful algebraic identities including the orthonormality relation and the Biedenharn-Elliott identity (see \cite[Appendix F]{TVi5}).

\subsection{Modular categories}\label{sect-modular}
Let $\bb$ be a ribbon fusion category and $I$ be a representative set of simple objects of $\bb$. The \emph{$S$-matrix} of~$\bb$ is the matrix $S=[S_{i,j}]_{i,j \in I}$, where
$$
\psfrag{i}[Br][Br]{\scalebox{.9}{$i$}}
\psfrag{j}[Bl][Bl]{\scalebox{.9}{$j$}}
S_{i,j}=\tr(c_{j,i}c_{i,j})=\left\langle  \rsdraw{.5}{.9}{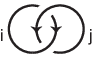} \right\rangle_\bb\in \End_\cc(\un)=\kk.
$$
Note that for any $i \in I$, the twist $\theta_i\co i \to i$ is multiplication by an invertible scalar
$v_i\in \kk$. We set
$$
\Delta_{\pm} =\sum_{i\in I} v_i^{\pm 1} \dim (i)^2\in \kk.
$$

A \emph{modular category} (over $\kk$) is a ribbon fusion category (over $\kk$) such that its $S$-matrix is invertible (over $\kk$).
If~$\bb$ is a modular category, then its dimension $\dim(\bb)$ and the scalars $\Delta_\pm$ are nonzero and satisfy $ \Delta_+\Delta_-=\dim (\bb)$, see~\cite{Tu1}.
We say that a modular category $\bb$ is {\it anomaly free} if
$\Delta_+=\Delta_-$.

Examples of modular categories are derived from quantum groups. The universal enveloping algebra $U\mathfrak{g}$ of a (finite dimensional complex) simple Lie algebra $\mathfrak{g}$ admits a deformation $U_q \mathfrak{g}$, which is a quasitriangular Hopf algebra. The representation category $\text{Rep}(U_q \mathfrak{g})$ is $\CC$-linear and ribbon. For generic $q \in \CC$, this category is semisimple. (The irreducible representations of $\mathfrak{g}$ can be deformed to irreducible representations of $U_q\mathfrak{g}$.) For $q$ an appropriate root of unity, a certain subquotient of $\text{Rep}(U_q\mathfrak{g})$ is a modular category with ground field $\kk = \CC$. For $\mathfrak{g} = \mathfrak{sl}_2(\CC)$, this result was pointed out by Reshetikhin and Turaev; the general case involves the theory of tilting modules.

Given a modular category $\bb$, the invariant $\langle \cdot \rangle_\bb$ of $\bb$\trait colored framed links and tangles extends by linearity to the case where colors are finite linear combinations of objects of $\bb$ with coefficients in $\kk$. In particular, the linear combination $$\Omega = \sum_{i \in I} \dim(i)\, i,$$  called the \emph{Kirby color}, has the following sliding property:
$$
\psfrag{U}[Bl][Bl]{\scalebox{.9}{$\Omega$}}
\psfrag{X}[Bl][Bl]{\scalebox{.9}{$X$}}
 \rsdraw{.5}{.9}{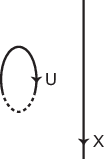} \,= \quad \rsdraw{.5}{.9}{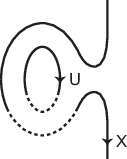}
$$
for any object $X$ of $\bb$  (meaning that the two tangles yield the same endomorphism of~$X$ under  $\langle \cdot \rangle_\bb$). Here, the dashed line represents an arc on the closed component colored by~$\Omega$. This arc can be knotted or linked to other components of the tangle (not shown in the figure). Also
$$
\langle O_\Omega^\pm \rangle_\bb = \Delta_\pm
$$
for the trivial knot $O^\pm_\Omega$ with framing $\pm 1$ and color~$\Omega$.

\section{Three dimensional TQFTs}\label{sect-TQFTs}
Inspired by the works of Witten \cite{Wi} and Segal \cite{Se},  Atiyah axiomatized in~\cite{At} the notion of a topological quantum field theory (TQFT). A 3-dimen\-sional TQFT~$Z$ (over a field $\kk$) assigns to every oriented closed surface $\Sigma$ a finite dimen\-sional $\kk$-vector space $Z(\Sigma)$ and assigns to every cobordism $(M, \Sigma, \Sigma')$ a $\kk$-linear homomorphism
$$
Z(M) = Z(M, \Sigma, \Sigma')\co Z(\Sigma) \to Z(\Sigma').
$$
Here, a \emph{cobordism} $(M, \Sigma, \Sigma')$ between two oriented closed surfaces $\Sigma$ and $\Sigma'$ is an oriented compact 3-manifold $M$ such that $\partial M = (-\Sigma) \sqcup \Sigma'$, where the boundary is oriented using the first outward pointing convention and the minus sign indicates the orientation reversal. A TQFT has to satisfy axioms which can be expressed by saying that
$$
Z\co \mathrm{Cob}_3 \to \mathrm{Vect}_\kk
$$
is a symmetric monoidal functor. Here $\mathrm{Vect}_\kk$ is the category of $\kk$-vector spaces and~$\mathrm{Cob}_3$ is the category whose objects are oriented closed surfaces, whose morphisms are diffeomorphism classes of cobordisms, and whose monoidal structure is given by the disjoint union.
In particular $Z(\emptyset)\cong\kk$ (where  $\emptyset$ is the empty surface) and
$$
Z( \Sigma  \sqcup \Sigma') \cong  Z(\Sigma ) \otimes Z( \Sigma')
$$ for any oriented closed
surfaces $\Sigma, \Sigma'$ (and similarly    for  cobordisms). Homeomorphisms of surfaces should induce isomorphisms of the corresponding vector spaces compatible with the action of cobordisms. Every oriented compact 3-manifold $M$ is a cobordism between $\emptyset$ and $\partial M$ so that~$Z$ yields a ``vacuum" vector $$Z(M) \in \Hom_\kk(Z(\emptyset), Z(\partial M))=Z(\partial M).$$ If $\partial M= \emptyset$, then this gives a numerical invariant $Z(M) \in \End_\kk\bigl(Z(\emptyset)\bigr)=\kk$.

An {\it isomorphism} of  3-dimensional TQFTs  $Z_1\to Z_2$  is a  natural monoidal isomorphism of   functors.
In particular, if two TQFTs $Z_1$, $Z_2$ are isomorphic, then $Z_1(M)= Z_2(M)$ for any  oriented closed 3-manifold $M$.

Interestingly, TQFTs are often defined for surfaces and 3-cobordisms with additional structure. The surfaces $\Sigma$ are normally endowed with Lagrangians, that is, with maximal isotropic subspaces in $H_1(\Sigma;\RR)$. For 3-cobordisms, several additional structures are considered in the literature: for example, 2-framings, $p_1$-structures, and numerical weights. All these choices are equivalent. The TQFTs requiring such additional structures are said to be \emph{projective} since they provide projective linear representations of the mapping class groups of surfaces, see \cite{Tu1}.

\section{The surgery approach}\label{sect-surgery-TQFT}
The Witten-Reshetikhin-Turaev invariants of oriented closed 3-manifolds are defined from modular categories and extend to 3-dimensional TQFTs. Their construction is based on the surgery presentation of 3-manifolds. In this section, we fix a modular category $\bb$ over a field $\kk$.

\subsection{Surgery on framed links}
Given an embedded solid torus $g \co S^1 \times D^2 \hookrightarrow S^3$,
where~$D^2$ is a 2-disk and $S^1 = \partial D^2$, a 3-manifold can be built as follows. Remove from $S^3$ the interior of $g(S^1 \times D^2)$ and glue back the solid torus $D^2 \times S^1$ along $g|_{S^1 \times S^1}$. This process is known as ``surgery". The resulting 3-manifold depends only on the isotopy class of the framed knot represented by $g$. More generally, surgery on a framed link $L = \cup_{i=1}^m L_i$ in $S^3$ with $m$ components yields an oriented closed 3-manifold $M_L$.

A theorem of Lickorish and Wallace asserts that any closed connected oriented 3-manifold is homeomorphic to $M_L$ for some $L$. Kirby proved that two framed links give rise to homeomorphic 3-manifolds if and only if these links are related by isotopy and a finite sequence of geometric transformations called \emph{Kirby moves}. There are two Kirby moves: adjoining a distant unknot $O^{\pm}$ with framing $\pm 1$ and sliding a link component over another one (as in the figure of the sliding property in Section~\ref{sect-modular}).

\subsection{The WRT invariants of closed 3-manifolds}
Let $L = \cup_{i=1}^m L_i$ be a framed link in $S^3$. Its linking matrix $(b_{i,j})_{1 \leq i,j \leq m}$ as coefficients defined as follows: for any~$i \neq j$, $b_{i,j}$ is the linking number of $L_i$ with $L_j$, and $b_{i,i}$ is the framing number of $L_i$. Denote by $e_+$ (resp. $e_{-}$) the number of positive (resp. negative) eigenvalues of this matrix. The sliding property of modular categories implies the following theorem. In its statement, a framed knot $K$ $\bb$-colored by the Kirby color $\Omega$ of $\bb$ is denoted by $K(\Omega)$.

\begin{thm}\label{thm-surgery-closed}
The expression
$$
\WRT_{\bb}(M_L)= \Delta_+^{-e_+} \Delta_-^{-e_{-}} \bigl\langle L_1(\Omega)\cup\dots \cup L_m(\Omega)  \bigr\rangle_\bb \in \kk
$$
is invariant under the Kirby moves on $L$. This expression yields, therefore, a well-defined topological invariant $\WRT_{\bb}$ of closed connected oriented 3-manifolds.
\end{thm}
Theorem~\ref{thm-surgery-closed} was first proved in \cite{RT} (see also \cite{Tu1}). In particular, the invariance under the second Kirby move follows from the sliding property of the Kirby color of a modular category (see Section~\ref{sect-modular}).
Several competing normalizations of $\WRT_{\bb}$ exist in the literature. Here, the normalization used is such that
$$\WRT_{\bb}(S^3) = 1 \quad \text{and} \quad \WRT_{\bb}(S^1 \times S^2) = \dim(\bb).$$ The invariant $\WRT_{\bb}$ extends to 3-manifolds  with a framed oriented $\bb$-colored link~$K$ inside (Wilson loops) by setting
$$
\WRT_{\bb}(M_L,K) = \Delta_+^{-e_+} \Delta_-^{-e_{-}} \bigl\langle L_1(\Omega) \cup \dots \cup L_m(\Omega)\cup K  \bigr\rangle.
$$

\subsection{The surgery TQFT}\label{sect-surg-TQFT-construction}
The Witten-Reshetikhin-Turaev invariants extend to a projective 3-dimensional TQFT denoted $\tau_{\bb}$ and called the \emph{surgery TQFT}. It depends on the choice of a square root~$\dd$ of $\dim(\bb)$. The \emph{coend} of the category $\bb$ is the object $$C=\bigoplus_{i \in I} i^* \otimes i,$$ where $I$ is a representative set $I$ of simple objects of $\bb$. For a connected oriented closed surface $\Sigma$ of genus $g$,
$$
\tau_{\bb}(\Sigma) =\Hom_{\bb}( \un, C^{\otimes g}).
$$
The dimension of this vector space enters the Verlinde formula
$$
\dim_{\kk}(\tau_{\bb}(\Sigma)) \, 1_{\kk}= 
\dim(\bb)^{g-1}
\sum_{i\in I}\dim(i)^{2-2g}
$$
where $1_\kk \in \kk$ is the unit of the field $\kk$. If $\text{char}(\kk) = 0$, then this formula computes $\dim_{\kk}(Z_{\bb}(\Sigma))$. For a closed connected oriented 3-manifold $M$ with numerical weight zero, $$\tau_{\bb}(M) = \dd^{-b_1(M)-1}\WRT_{\bb}(M),$$ where $b_1(M)$ is the first Betti number of $M$. In particular,
$$
\tau_{\bb}(S^3) = \dd^{-1} \quad \text{and} \quad \tau_{\bb}(S^1 \times S^2) = 1.
$$
The two dimensional part of $\tau_{\bb}$ determines a ``modular functor" in the sense of Segal, Moore, and Seiberg.

The TQFT $\tau_{\bb}$ extends to a vaster class of surfaces and cobordisms. Surfaces may be enriched with a finite set of marked points, each colored with an object of $\bb$ and endowed with a tangent direction. Cobordisms may be enriched with ribbon (or fat) graphs whose edges are colored with objects of $\bb$ and whose vertices are labeled with appropriate intertwiners.  The resulting TQFT, called the \emph{surgery graph TQFT} and also denoted $\tau_{\bb}$, is nondegenerate in the sense that, for any surface $\Sigma$, the vector space~$\tau_{\bb}(\Sigma)$ is spanned by
the vacuum vectors determined by all~$M$ with $\partial M = \Sigma$.  A detailed construction of~$\tau_{\bb}$ is given in \cite{Tu1}.

If $\bb$ is anomaly free and $\dd=\Delta_\pm$, then $\tau_{\bb}$  is a genuine 3-dimensional TQFT (not only a projective one).

\section{The state sum approach}\label{sect-state-sum-TQFT}
Another approach to three dimensional TQFTs is based on the theory of $6j$-symbols and state sums on triangulations of 3-manifolds. This approach, introduced by Turaev and Viro in 1992 and refined by Barrett-Westbury in 1995, is a quantum deformation of the Ponzano-Regge model for three dimensional lattice gravity. The state sum quantum invariants of closed 3-manifolds are defined from spherical fusion categories with nonzero dimensions and extend to 3-dimensional TQFTs. 

In this section, we fix a spherical fusion category $\cc$ (over a field $\kk$) with nonzero dimension and let $I$ be a representative set of simple objects of  $\cc$.

\subsection{Triangulations of 3-manifolds}
A \emph{tetrahedron} is the convex hull of four affinely independent points in some affine space. It has 4 triangular faces called \emph{triangles}, 6 edges, and 4 vertices:
$$
\rsdraw{.45}{.9}{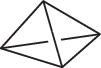} \;.
$$

A \emph{triangulation} of a 3-manifold $M$ is a decomposition of~$M$ into finitely many tetrahedra such that the triangles of the tetrahedra are identified with each other pairwise, and the interiors of the tetrahedra remain disjoint.

Moise proved that any compact 3-manifold has a triangulation. Pachner proved that two triangulations of a 3-manifold are related by a finite sequence of ambient isotopies of triangulations, 2-3 moves, 1-4 moves, and their inverses.
The 2-3 move is performed on two different tetrahedra meeting in a
triangle. It deletes this triangle by introducing a new edge connecting the opposite
corners of the tetrahedra (creating  three new tetrahedra):
$$
\rsdraw{.45}{.9}{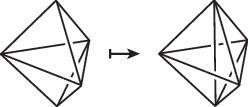} \;\,.
$$
The 1-4 move introduces a vertex inside a tetrahedron and connects it to the four vertices of the tetrahedron with four edges (creating  four new tetrahedra):
$$
\rsdraw{.45}{.9}{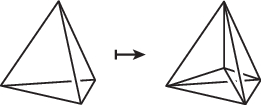} \;\,.
$$

\subsection{State sum invariants of closed 3-manifolds}\label{set-state-sum-closed}
Let $M$ be an oriented closed 3-manifold.
Pick a  triangulation  of~$M$ and a total order on the set of vertices of the triangulation. A \emph{state} is a map from the set of edges of the triangulation to $I$. Note that the number of states is finite since both the set of edges and $I$ are finite.
For a given state $s$, we set
$$
\dim(s)=   \prod_{e} \dim (s(e))    \in \kk
$$
where $e$ runs over all edges of the triangulation.
Next, we define    a    scalar   $\vert s \vert \in \kk$ as follows.

For any triangle $t$ of the triangulation, consider the $\kk$-vector spaces
\begin{align*}
& H_{s,t}^+=\Hom_\cc\bigl(s(02),s(01) \otimes s(12)\bigr),\\
& H_{s,t}^-=\Hom_\cc\bigl(s(01) \otimes s(12),s(02)\bigr),
\end{align*}
where $0<1<2$ are the vertices of $t$ and $(ij)$ denotes the edge connecting the vertices $i$ and $j$. Since the category $\cc$ is fusion, the pairing $$\alpha \otimes \beta \in H_{s,t}^- \otimes H_{s,t}^+ \mapsto \tr(\alpha\beta) \in \kk$$  is non-degenerate. We denote by $\ast_{s,t}$ the image of $1_\kk$ under its inverse copairing $\kk \to  H_{s,t}^+ \otimes H_{s,t}^-$. Let
$$
H_s=   \underset{t}{\bigotimes} \, H_{s,t}^+ \otimes H_{s,t}^-
$$
be the unordered tensor product of $H_{s,t}^+$ and $H_{s,t}^-$ over all triangles  $t$  of the triangulation, and set
$$
\ast_s=\bigotimes_t \ast_{s,t} \in  H_s.
$$

For any tetrahedron $\Delta$ of the triangulation, set $\varepsilon_\Delta=+$  if the orientation of~$\Delta$ induced by the order of its vertices coincides with that induced by $M$,  and set $\varepsilon_\Delta=-$ other\-wise. Moreover, given any triangle $t$ at the boundary of~$\Delta$, set $\varepsilon(t,\Delta)=+$ if the orientation of $t$ induced by the order of its vertices coincides with the boundary orientation of $t \subset \partial \Delta$ induced by the orientation of $M$ restricted to~$\Delta$, and set $\varepsilon(t,\Delta)=-$ otherwise. Section~\ref{sect-6j} yields the $6j$\trait symbol
$$
|\Delta|_s=\begin{Bmatrix} s(01) & s(02) & s(12) \\ s(23) & s(13) & s(03) \end{Bmatrix}_{\!\varepsilon_\Delta} \co  \bigotimes_{t \subset \partial \Delta}  H_{s,t}^{\varepsilon(t,\Delta)} \to \kk
$$
where $t$ runs over all triangles in the boundary of $\Delta$.

Since $M$ is closed, each triangle $t$ of the triangulation is adjacent to two tetrahedra $\Delta_1$ and~$\Delta_2$ of the triangulation and $\varepsilon(t,\Delta_2)=-\varepsilon(t,\Delta_1)$. Then the unordered tensor product over all tetrahedra~$\Delta$ of their associated $6j$-symbols is a $\kk$-linear form
$$
V_s= \bigotimes_\Delta  |\Delta|_s  \co H_s \to \kk.
$$
Evaluating  $V_s$ on   $\ast_s$ yields $ \vert s\vert=  V_s(\ast_s) \in \kk$. Finally,   set
$$
|M|_\cc=\dim (\cc)^{-\upsilon} \, \sum_{s} \,  \dim (s) \,   \vert s\vert  \in \kk,
$$
where $s$ runs over all states of the triangulation of $M$ and~$\upsilon$  is the number of vertices of the triangulation.

\begin{thm}\label{thm-state-3man}
$|M|_\cc$ is a topological invariant of $M$  independent of the
choice    of~the triangulation and~$I$.
\end{thm}
For example, one computes that
$$
|S^3|_\cc=\dim(\cc)^{-1} \quad \text{and} \quad |S^1 \times S^2|_\cc=1.
$$

When $\cc$ is the fusion category derived from the representations of the quantum group $U_q(\mathfrak{sl}_2\CC)$ with $q$ an appropriate root of unity (see Section~\ref{sect-modular}), then $|M|_\cc$ is equal to the original Turaev-Viro invariant \cite{TV} of $M$.

The proof of Theorem~\ref{thm-state-3man}  consists in particular of verifying the invariance of the state sum under the application of Pachner moves on the triangulation. This comes down to the orthonormality relation and the Biedenharn-Elliott identity for $6j$-symbols, see  \cite{BW}.

The state sum may be more generally defined on skeletons of 3-manifolds (including triangulations, their dual cellular decompositions, and spines), see \cite{TVi5}.

\subsection{The state sum TQFT}
If $M$ is an oriented compact 3-manifold with nonempty boundary, then the algorithm described in the previous section applied to a state $s$ of a triangulation of $M$ yields not a scalar but a $\kk$-linear form
$$
\vert s\vert \co H_s^{\partial}=   \underset{t}{\bigotimes} \, H_{s,t}^{\varepsilon(t,\Delta_t)} \to \kk
$$
where  $t$  runs over all triangles in  the boundary of $M$ and~$\Delta_t$ denotes the unique tetrahedron adjacent to such a triangle~$t$.
Consider the state sum
$$
|M|^\circ=\dim (\cc)^{-\upsilon} \, \sum_{s} \,  \dim (s) \,   \vert s\vert
$$
where $s$ runs over all states and~$\upsilon$  is the number of vertices in the interior of $M$.
Then the assignment $M \mapsto |M|^\circ$ behaves well with the gluing of 3-manifolds along boundary components.
Consequently, there is a standard procedure (see \cite{TVi5}) to transform it into a (genuine) TQFT
$$
\vert \cdot \vert_\cc \co \mathrm{Cob}_3 \to \mathrm{Vect}_\kk.
$$

For example, the vector space associated to the 2-sphere is $\vert S^2 \vert_\cc \cong \kk$.

\section{Comparison of the two approaches}\label{sect-comparison}

The comparison of the surgery and state sum appro\-aches to 3-dimensional TQFTs is based on the notion of center of a monoidal category due to Joyal,  Street, and Drinfeld.

\subsection{Categorical centers}
The \emph{center} of a monoidal category $\cc$ is the braided category $\zz(\cc)$ defined as follows. The objects of~$\zz(\cc)$ are \emph{half braidings} of $\cc$, that is, pairs $(A,\sigma)$, where $A$ is an object of $\cc$ and
$
\sigma=\{\sigma_X \co  A \otimes X\to X \otimes A\}_{X \in \cc}
$
is a natural isomorphism such that
$$
 \sigma_{X \otimes Y}=(\id_X \otimes
\sigma_Y)(\sigma_X \otimes \id_Y).
$$
A morphism $(A,\sigma)\to (A',\sigma')$ in $\zz(\cc)$ is a
morphism $f \co A \to A'$ in $\cc$ such that $$(\id_X \otimes f)\sigma_X=\sigma'_X(f \otimes \id_X)$$ for all $X\in \cc$. 
The
  unit object of $\zz(\cc)$ is $\un_{\zz(\cc)}=(\un,\{\id_X\}_{X \in
\cc})$ and the monoidal product is
\begin{equation*}
(A,\sigma) \otimes (B, \rho)=\bigl(A \otimes B,(\sigma \otimes \id_B)(\id_A \otimes \rho) \bigr).
\end{equation*}
The braiding $c$ in $\zz(\cc)$ is defined by
$$
c_{(A,\sigma),(B, \rho)}=\sigma_{B}\co (A,\sigma) \otimes
(B, \rho) \to (B, \rho) \otimes (A,\sigma).
$$

If $\cc$ is a $\kk$-linear category, then so is $\zz(\cc)$. If $\cc$ is pivotal, then so is $\zz(\cc)$ with $(A,\sigma)^*=(A^*,\sigma^\dag)$, where
$$
 \psfrag{M}[Bc][Bc]{\scalebox{.9}{$A$}}
 \psfrag{X}[Bc][Bc]{\scalebox{.9}{$X$}}
 \psfrag{s}[Bc][Bc]{\scalebox{.9}{$\sigma_{X^*}$}}
\sigma^\dag_X= \rsdraw{.45}{1}{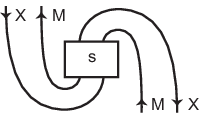} \co A^* \otimes X \to X \otimes A^*,
$$
and (co)evaluations morphisms and pivotal structure are inherited from $\cc$.
The (left and right) traces of morphisms and   dimensions of objects
in~$\zz(\cc)$ are the same as in~$\cc$.

\subsection{The comparison}
The first connections between the surgery and state sum constructions were established by Walker \cite{Wa} and Turaev \cite{Tu1}: if $\bb$ is a modular category, then it is also a spherical category with nonzero dimension and the surgery and state sum invariants
are related by:
\begin{equation}\label{eq-main-mod}
\vert
M\vert_{\bb}=\tau_{\bb} (M)\,\tau_{\bb} (-M)
\end{equation}
for every oriented closed 3-manifold $M$, where $-M$ is the $3$-manifold $M$ with opposite orientation. In particular, if~$\bb$ is unitary over $\kk=\CC$ (meaning that the Hom spaces in $\bb$ are equipped with a conjugation compatible with the pivotal structure and the braiding), then $\tau_{\bb} (-M)=\overline{\tau_{\bb} (M)}$ and so $\vert M\vert_{\bb}=|\tau_{\bb} (M)|^2 \in \RR_+$.

But in general, a spherical category need not be braided and so cannot be used as input to define the Witten-Reshetikhin-Turaev invariant.
However, let $\cc$ be   a  spherical fusion  category over an algebraically closed field~$\kk$  such that $\dim (\cc)\neq 0$. A fundamental theorem of M\"uger~\cite{Mu} asserts that   the center $\zz(\cc)$ of $\cc$ is an anomaly free modular category with $\Delta_+=\Delta_-=\dim (\cc)$. In particular
$$
\dim \bigl(\zz(\cc)\bigr)=\Delta_+ \Delta_-= \dim (\cc)^2.
$$
Consequently, such a category $\cc$  gives rise to two (genuine) 3-dimen\-sio\-nal TQFTs: the state sum TQFT $\vert \cdot \vert_{\cc}$ and the surgery TQFT~$\tau_{Z(\cc)}$ associated with the square root $\dim (\cc)$ of $\dim(\zz(\cc))$.

\begin{thm}\label{thm-main-closed}
The TQFTs $\vert \cdot \vert_\cc$ and $\tau_{\zz(\cc)}$ are isomorphic. In particular, for any  oriented closed 3-manifold~$M$,
\begin{equation}\label{eq-main}
\vert M\vert_{\cc}=\tau_{Z(\cc)} (M),
\end{equation}
and for any oriented closed surface $\Sigma$,
\begin{equation}\label{eq-main2}
\vert \Sigma\vert_{\cc}\cong \tau_{Z(\cc)} (\Sigma).
\end{equation}
\end{thm}

Theorem~\ref{thm-main-closed}  was first proved in~\cite{TVi0} (see also \cite{TVi5}). In the case where the characteristic of~$\kk$ is equal to zero, Theorem~\ref{thm-main-closed} was independently proved  in \cite{Ba2}.

Theorem~\ref{thm-main-closed} relates through the categorical center two categorical approaches to invariants of 3-manifolds. This relationship sheds new light on both approaches and shows, in particular, that the  surgery approach is more general than the state sum approach. Formula \eqref{eq-main2} gives
$$
\vert \Sigma\vert_{\cc}\cong \Hom_{\zz(\cc)}\bigl(\un_{\zz(\cc)},\mathbb{C}^{\otimes g}\bigr)
$$
where $\mathbb{C}$ is the coend of $\zz(\cc)$ and $g$ is the genus of $\Sigma$. Note that $\mathbb{C}=(A,\sigma)$ can be computed explicitly using the category $\cc$, see \cite{TVi5}. In particular, $$A=\bigoplus_{i,j\in I} i^* \otimes j^* \otimes i \otimes j.$$

The formula \eqref{eq-main} was previously known to be true in several special cases:  when~$\cc$ is the category of representations of a finite group, when~$\cc$ is the category of bimodules associated with a subfactor \cite{KSW}, and when~$\cc$   is modular \cite{Tu1,Wa}. In the latter case,  Formula \eqref{eq-main-mod} can indeed be derived from Formula \eqref{eq-main}: if $\bb$ is a modular category, then its center $\zz(\bb)$ is braided equivalent to the Deligne tensor product $\bb\boxtimes \overline{\bb}$ (where $\overline{\bb}$ is the mirror of $\bb$) and therefore Formula~\eqref{eq-main} gives $$\vert M\vert_{\bb}=\tau_{\bb \boxtimes \overline{\bb}} (M)=\tau_{\bb} (M) \,\tau_{\overline{\bb}} (M)=\tau_{\bb} (M)\,\tau_{\bb} (-M).$$

\section{Generalizations and perspectives}\label{sect-generalizations}
Three dimensional TQFTs have several interesting generalizations including extended TQFTs, defect TQFTs, homotopy QFTs, and non-semisimple TQFTs. Extended TQFTs are motivated by applications of higher categorical ideas to the functorial (cut-paste) nature of TQFTs while defect TQFTs incorporate the presence of defects in the underlying manifolds, which originate from certain concepts in physics such as domain walls, boundaries, and interfaces. On the other hand, homotopy QFTs can be seen as TQFTs for manifolds endowed with an extra structure encoded by a homotopy class of maps to a target space (viewed as the classifying space of the structure).
Non-semisimple theories weaken the semisimplicity condition on the underlying fusion categories and overcome certain limitations on the quantum invariants.

In the following subsections, we shortly discuss these generalizations and review some recent works in these fields.

\subsection{Extended TQFTs}
Recall that a $3$-dimensional TQFT provides a numerical invariant of oriented closed $3$-mani\-folds. This invariant can be computed by cutting the $3$-mani\-fold along codimension one submanifolds into $3$-manifolds with boundary, and then by composing the corresponding linear maps.
However, these maps as well as the vector spaces assigned to the boundary surfaces, are not always easy to determine. This leads to wanting to cut the $3$-manifold along higher codimensional submanifolds. This motivates the definition and study of extended TQFTs.

While a $3$-dimensional TQFT assigns algebraic invariants to closed surfaces and compact $3$-manifolds, a once extended $3$-dimensional TQFT should assign  algebraic invariants to  closed $1$-manifolds, $2$-mani\-folds, and  $3$-mani\-folds with corners. More precisely, there is a symmetric monoidal $2$-category $\mathrm{Bord}_{3,2,1}$ which extends the category $\mathrm{Cob}_3$ and
whose objects are oriented closed $1$-manifolds, $1$\trait morphisms are oriented 2-dimensional cobordisms between them, and $2$-morphisms are diffeomorphism classes of oriented cobordisms between $1$-morphisms (such cobordisms are oriented compact $3$-manifolds with codimension~$2$ corners).
A \emph{once extended} $3$-dimensional TQFT is then a symmetric monoidal  2-functor  from $\mathrm{Bord}_{3,2,1}$ to some algebraic symmetric monoidal $2$-category. For example, the Witten-Reshetikhin-Turaev surgery graph TQFT (see Section~\ref{sect-surg-TQFT-construction}) can be seen as a once extended $3$-dimen\-sional TQFT (with anomaly).

A classification of once extended $3$-dimensional TQFTs with values in the symmetric monoidal 2-category of Cau\-chy complete
linear categories (over an  algebraically closed field) is given in \cite{BDSPV2}.
This classification states a one-to-one correspondence between equivalence classes of such extended TQFTs and equivalence classes of modular tensor categories whose ano\-maly factor is $1$. This correspondence takes an extended TQFT to its value on the circle.

One can further try to extend a once extended $3$-dimen\-sional TQFT to three codimensional submanifolds, that is, to points. Such an extended TQFT is called \emph{fully extended}
and is formally defined as a symmetric monoidal $3$-functor from the $3$-category $\mathrm{Bord}_{3,2,1,0}$ to some algebraic symmetric monoidal $3$-category.
The cobordism hypothesis, conjectured by Baez-Dolan \cite{BD} and proven by Lurie \cite{Lu} and recently by Ayala-Francis \cite{AF}, states that a fully extended framed TQFT is determined by its value on a point, and any fully dualizable object of the target category gives rise to a fully extended framed TQFT which assigns that object to a point, see \cite{Fr,DSS}. Moreover, Lurie \cite{Lu} generalized the cobordism hypothesis to arbitrary tangential structures on manifolds by using homotopy fixed points. In this formulation, fully extended oriented $3$-dimensional TQFTs are classified by homotopy $SO(3)$-fixed points of the target 3-category.

A natural candidate for the target symmetric monoidal $3$-category is
the $3$-category $\mathrm{TC}$ whose objects are finite rigid monoidal linear categories, $1$-morphisms are finite bimodule categories, $2$-morphisms are bimodule functors, and $3$-morphisms are bimodule natural transformations. In this case, fully dualizable objects and homotopy $SO(3)$-fixed points in $\mathrm{TC}$ are computed in \cite{DSS} as fusion categories of nonzero dimension and spherical fusion categories, respectively. Given a spherical fusion category $\cc$ of nonzero dimension, the associated fully extended oriented $3$-dimensional TQFT conjecturally extends the state sum  TQFT $\vert \cdot \vert_\cc$ associated with~$\cc$ (see Section~\ref{sect-state-sum-TQFT}).

\subsection{TQFTs with defects}
Defect TQFTs generalize TQFTs by allowing the presence of defects, which are lower dimensional submanifolds of the cobordisms that can carry nontrivial topological or quantum information. More precisely, a \emph{$3$-dimensional defect TQFT} is a symmetric monoidal functor
$$
Z \co \mathrm{Cob}^{\text{def}}_3 (\mathbb{D}) \to \Vect_{\kk},
$$
where $\mathrm{Cob}^{\text{def}}_3(\mathbb{D})$ is the category of oriented closed surfaces and oriented cobordisms endowed with a stratification by submanifolds labeled with elements of a fixed labeling data~$\mathbb{D}$.
The labelings of defect submanifolds should satisfy (higher categorical) algebraic relations reflecting their adjacency.
In particular, it is shown in \cite{CMS} that any $3$-dimensional defect TQFT yields a $\kk$-linear Gray 3\trait cate\-gory with duals.

Concrete examples of $3$-dimensional defect TQFTs generalizing the surgery and state sum TQFTs
have been developed in various papers, notably \cite{KS1,FSV,KK,CRS1,Me}.

\subsection{Homotopy QFTs}
Roughly, HQFTs are TQFTs for manifolds endowed with a map to a fixed target space. More precisely, let $X$ be a connected topological space. Following Turaev \cite{Tu}, a $3$-dimensional \emph{homotopy quantum field theory} (HQFT) \emph{with target~$X$} is a symmetric monoidal functor $$Z \co X\mathrm{Cob}_3 \to \mathrm{Vect}_\kk.$$ Here $X\mathrm{Cob}_3$ is the symmetric monoidal category whose objects are oriented closed surfaces endowed with a map to $X$ and whose morphisms are diffeomorphism classes of oriented cobordisms equipped with a homotopy class of maps to $X$ restricting to the given maps on their boundary.  In particular, a $3$-dimensional HQFT with target~$X$ produces a scalar homotopy invariant of maps from oriented closed $3$-manifolds to $X$. For example, any  third cohomology class $\theta \in H^3(X,\kk^*)$ gives rise to a 3-dimensional HQFT with target $X$, called \emph{cohomological HQFT},  whose associated homotopy invariant of a map $f\co M \to X$ is the evaluation of the pullback class $f^*(\theta) \in H^3(M,\kk^*)$ with the fundamental class $[M] \in H_3(M,\Z)$.

If $X$ is a connected homotopy 0-type (that is, a contractible space), then any HQFT with target $X$ is equivalent to a TQFT.

If $X$ is a connected homotopy 1-type, then $X$ is a $K(G,1)$ space where $G$ is the fundamental group of $X$. In this case, the surgery and state sum TQFTs have been generalized in \cite{TVi1,TVi3} to 3-dimensional HQFTs with target $X$. The relevant algebraic structures for their construction are modular and spherical fusion categories which are $G$-graded (meaning that objects have a degree in $G$ and this degree is multiplicative with respect to the monoidal product). Generalizing Theorem~\ref{thm-main-closed}, it is shown in \cite{TVi6} that the surgery and state sum HQFTs are  related via the graded center of graded fusion categories. Also, the orbifold construction \cite{SW2} associates a TQFT to each HQFT with target $X$. For example,  the Dijkgraaf-Witten TQFT is the orbifoldization of a cohomological HQFT.

If $X$ is a  connected  homotopy 2-type, then $X$ may be encoded by a crossed module $\chi \co E \to H$  which is a certain group homomorphism with $\pi_1(X)=\text{coker}(\chi)$ and $\pi_2(X)=\ker(\chi)$. In this case, the state sum TQFT has been generalized  in \cite{SV}  to a 3-dimensional HQFT with target~$X$. For this purpose, the relevant algebraic inputs are $\chi$-graded spherical fusion categories. These are a class of monoidal categories in which not only the objects have a degree (in $H$) but also the morphisms have a degree (in~$E$), and the compatibility of these degrees is governed by the crossed module $\chi$. For example, the cohomological HQFTs associated with~$X$ are particular instances of state sum  HQFTs with target~$X$.

\subsection{Non-semisimple quantum invariants}
Hennings \cite{He} was the first to build a non-semisimple quantum invariant  of closed $3$-manifolds by using a finite dimensional ribbon Hopf algebra. When the Hopf algebra is semisimple, this invariant agrees with  the  Witten-Reshetikhin-Turaev invariant derived from the category of representations of the Hopf algebra. Lyubashenko \cite{Ly} extended Hennings' construction by using ribbon finite tensor categories. Note that the Lyubashenko  invariant   does not form a TQFT in the usual sense because it does not behave well under the disjoint union operation (in particular, when the category is not semisimple, it vanishes on  all closed $3$-manifolds with positive first Betti number). However, the Lyubashenko  invariant forms an extended TQFT in a weaker sense (by considering cobordisms with corners between connected surfaces and using the connected sum as monoidal product), see~\cite{KL}.

To construct genuine TQFTs from non-semisimple modular categories, a useful tool is that of a modified trace introduced in~\cite{GPV}.
Such traces have been used in \cite{CGP} to define a non-semisimple version of the surgery quantum invariants.
The CGP invariants are actually part of an extended TQFT for admissible cobordisms decorated with colored ribbon graphs and cohomology classes, see \cite{DR}.

Another instance of a non-semisimple invariant is the Kuperberg invariant \cite{Ku2} of framed 3-manifolds defined from any finite dimensional Hopf algebra by using Heegaard splittings (i.e.,  decompositions of 3-manifolds into two handlebodies). If the Hopf algebra is semisimple, then the Kuperberg invariant is an invariant of closed 3\trait manifolds and is equal (by Theorem~\ref{thm-main-closed}) to the Hennings invariant derived from the Drinfeld double of the Hopf algebra. This result is extended to non-semisimple
Hopf algebras in \cite{CC}.

A non-semisimple generalization of the Turaev-Viro state sum invariant of closed 3-manifolds is given in \cite{CGPT} using a spherical finite tensor category as algebraic input. It is extended to a (non-compact) 3-dimensional TQFT in  \cite{CGPV} via Juh\'asz's presentation  \cite{Ju} of the category $\mathrm{Cob}_3$ by generators and relations.

\subsection*{Acknowledgements}
This work was supported in part by the Labex CEMPI  (ANR-11-LABX-0007-01).

\end{document}